\documentclass[11pt, twoside]{amsart}
\usepackage[top=4cm, bottom=3cm, left=3cm, right=3cm, a4paper]{geometry}
\usepackage{graphicx,amssymb,amsrefs,amsmath,amsthm} 

\newcommand{\ud}[1]{\underline{#1}}
\newcommand{\Hom}{\mathrm{Hom}}
\newcommand{\sHom}{\ud{\Hom}}
\newcommand{\Mod}{\mathrm{Mod}}
\newcommand{\cA}{\mathcal{A}}
\newcommand{\scA}{\underline{\mathcal{A}}}
\newcommand{\Mor}{\mathrm{Mor}}
\newcommand{\Ker}{\mathrm{Ker}}
\newcommand{\sMor}{\ud{\Mor}}
\newcommand{\cC}{\mathcal{C}}
\newcommand{\cF}{\mathcal{F}}
\newcommand{\cP}{\mathcal{P}}
\newcommand{\Proj}{\mathrm{Proj}}
\newcommand{\Obj}{\mathrm{Obj}}

\newtheorem{theorem}{Theorem}[section]
\newtheorem{lemma}[theorem]{Lemma}
\newtheorem{corollary}[theorem]{Corollary}

\title{Weak Generators in Stable Module Categories}
\author{Christian Lomp}
\address{Department of Mathematics, University of Porto, 4169-007 Porto, Portugal}
\email{{clomp@fc.up.pt}}
\date{November 2026}

\begin{document}
\begin{abstract}
The aim of this note is to show that the existence of weak generators in the stable module category of a ring, or equivalently, of $p$-indigent modules in the terminology of Holston et al., depends on the underlying model of set theory.
\end{abstract}

\thanks{ I acknowledge partial support by CMUP, member of LASI, which is financed by national funds through FCT - Fundação para a Ciência e a Tecnologia, I.P., under the projects with reference UIDB/00144/2020 and UIDP/00144/2020.  Furthermore, I would like to thank Yilmaz Dur\u{g}un for stimulating discussions on subprojectivity during his visit to the University of Porto and for drawing my attention to the problems raised in \cite{Holston}. I am also grateful to Driss Bennis for insightful discussions on the categorical aspects of the subject. Finally, I acknowledge the use of AI for literature review which lead me to discover the excellent book by Eklof and Mekler, whose distinctive blue cover I had noticed many times in our library without ever having ventured to open it.
This paper is dedicated to Gerhard Lomp.}
\maketitle

\section{Introduction}
A generator in an additive category $\cC$ is an object $G$ such that $\Mor_{\cC}(G,-)$ is faithful, by which one understands that for any distinct morphisms  $f,g \in \Mor_{\cC}(X,Y)$ there exists $h\in \Mor_{\cC}(G,X)$, such that $f\circ h \neq g\circ h$. In this note  will call an object $G$ a weak generator if  $\Mor_{\cC}(G,X)$ is non-zero, for any non-zero $X\in \Obj(\cC)$. We are interested in the existence of weak generators in the stable category of an Abelian category. The stable category of an Abelian category $\cA$ with enough projectives is the quotient category $\scA$, where $\Obj(\scA) = \Obj(\cA)$ and $\sMor(X,Y)= \Mor(X,Y)/\sim$, with $f\sim g$ if and only if $f-g$ factors through some projective object $P\in \Proj(\cA)$. 
The objects of $\cA$ that are mapped to the zero objects in $\scA$ under the natural projection functor $F:\cA \to \scA$  are precisely the projective objects  $\Proj(\cA)$. A weak generator in $\scA$ is therefore an object $G$ such that for any non-projective object $X\in \cA$, there exists a morphism $f:G\to X$ that does not factor through a projective object. Our interest in the existence stems from an open problem in homological algebra asking whether \emph{p-indigent} modules exist over any (associative, unital) ring. More precisely, Holston et al defined in \cite{Holston} what they called the \emph{subprojectivity domain} of a (left) $R$-module $M$, in order to analyze to what extent $M$ is projective. A module $X$ belongs to that subprojectivity domain, denoted by $\mathfrak{Pr}^{-1}(M)$,  if for any epimorphism $f:Y\to X$ and homomorphism $g:M\to X$ can be lifted to a homomorphism $\overline{g}:M\to Y$, such that $g=f\overline{g}$.  Since any $R$-module is an epimorphic image of a free module, it is not difficult to see that $X\in \mathfrak{Pr}^{-1}(M)$ if and only if every $g\in \Hom(M,X)$ factors through a projective module (see \cite[Lemma 2.3]{Holston} or \cite[Lemma IV.1.5]{ARS}). Thus, $$\mathfrak{Pr}^{-1}(M) = \Ker(\sHom(M,-),$$
where by $\Ker(\sHom(M,-))$ we mean the class of $R$-modules $X$ with $\sHom(M,X)=0$.
Holston et al. defined a module $G$ to be $p$-indigent if  $\mathfrak{Pr}^{-1}(G)=Proj(R)$ or in other words, if $\sHom(G,-)$ is non-zero on non-projective modules, which simply means that $G$ is a weak generator in the stable category $ \ud{R\text{-}\Mod}$.
The question of whether  $p$-indigent modules exist   has not yet been fully resolved. A result by Yilmaz Dur\u{g}un says that if $\mathcal{FP}$ denotes a set of representatives of finitely presented left $R$-modules  and $G=\bigoplus_{X\in \mathcal{FP}} X$, then $\mathfrak{Pr}^{-1}(G)$ is equal to the class of flat modules. Hence, if $R$ is left perfect, or equivalently if all flat modules are projective, then $G$ is $p$-indigent (see \cite[Proposition 2.1]{yilmaz}). However, beyond this result, not much is known about the existence of $p$-indigent modules, or equivalently, weak generators of the stable module category. This note tries to show that the existence of such objects depends on the model of set theory we are choosing.

\section{Almost Projective Modules}
Let us reformulate the problem using the notion of almost projective modules.  For an infinite cardinal $\kappa$ we say that a module $M$ is $\kappa$-projective, if every submodule $N$ of $M$, with $|N|<\kappa$ is contained in a projective submodule of $M$.\footnote{This definition is not standard. For example, Braun and Trlifaj defined in \cite{BraunTrlifaj} a module   $M$ to be $\kappa$-projective, for a regular uncountable cardinal $\kappa$, if $M$ has a $\kappa$-dense system consisting of $<\kappa$-generated projective modules. Here, a direct system, $\cF$, of submodules of $M$ is said to be a $\kappa$-dense system in $M$ if $\cF$ is closed under unions of well-ordered ascending chains of length $< \kappa$, and every subset of $M$ of cardinality $< \kappa$ is contained in an element of $\cF$. Since in their definition any submodule of cardinality $<\kappa$ is contained in a projective module, these modules are $\kappa$-projective in the sense of this paper.}

\begin{lemma}\label{necessary} If the stable category of $R$ has a weak generator of cardinality $\lambda$, then any $\kappa$-projective module is projective, for every infinite cardinal $\kappa > \lambda$.
\end{lemma}
\begin{proof} Suppose $G$ is weak generator in $\ud{R\text{-}\Mod}$. Let $\lambda = |G|$ and $\kappa > \lambda$ an infinite cardinal. Let $M$ be any  $\kappa$-projective module $M$. Then for any  homomorphism  $f:G\to M$, we have that $|f(G)|<\kappa$. Since $M$ is $\kappa$-projective, $f(G)$ is contained in a projective submodule $P\subseteq M$ and $f$ factors through $P$. This means  $f=0$ in $\ud{R\text{-}\Mod}$ and $\sHom(G,M)=0$. As $G$ is a weak generator, $M$ must be projective.
\end{proof}

The Lemma could also be stated as follows: if $G$ is a weak generator of the stable category of left $R$-modules, then any $\max(|G|^+,\aleph_0)$-projective module is projective.

\begin{lemma}\label{sufficient} Let $R$ be a left hereditary ring. The stable category of $R$ has a weak generator if and only if there exists an infinite cardinal $\kappa$ such that every  $\kappa$-projective module is projective.
\end{lemma}
\begin{proof} Suppose that there exists an infinite cardinal $\kappa$ such that every  $\kappa$-projective module is projective.
Let $\cF$ be the set of representatives of non-projective modules $H$ with $|H|<\kappa$ and set $G=\bigoplus_{H\in \cF} H$. Let $M$ be any $R$-module with $\sHom(G,M)=0$ and take a submodule  $N\subseteq M$ of cardinality $|N|<\kappa$. If $N$ is not projective, then there exists an isomorphism $\alpha: H\to N$, for some $H\in \cF$.
Since $H$ is a direct summand of $G$, there is a projection $\pi: G\to H$ and an embedding  $\epsilon:H\to G$ satisfying $\pi\epsilon=id_H$. Thus, $f = \alpha \pi:G\to M$ is a homomorphism, which by assumption factors through a projective module $P$. Hence  there are homomorphisms $g:G\to P$ and $h:P\to M$ such that $f=hg$. Then $\alpha = f\epsilon = h (g\epsilon)$ and $g\epsilon: H\to P$ is an embedding. As $R$ is hereditary and $P$ projective, $H$ and $N$ would be projective, which is a contradiction to our assumption that $N$ is not. Thus, any submodule $N$ of $M$, with $|N|<\kappa$ is projective and therefore, $M$ is $\kappa$-projective and hence projective by assumption. This means that $G$ is a weak generator in $\ud{R\text{-}\Mod}$.
The converse holds by Lemma \ref{necessary}.
\end{proof}

Using well-known set theoretical results from Eklof and Mekler's book \cite{Eklof}, we show that by adding axioms to ZFC, the question of the existence of a weak generator in the stable category can be answered.\footnote{The reader might be warned that the author is not an expert in set theory.}

\section{Axiom of Constructibility}

In set theory the Axiom of Constructibility, for short $V=L$,  says that all sets are constructible. This axiom is  independent of ZFC, assuming ZFC is consistent. This means that neither can $V=L$ nor $V\neq L$ be proven within ZFC. For the reader's sake we recall the definitions of $V$ and $L$ (see \cite[VI.3]{Eklof}).The (von Neumann) universe $V$ of sets is defined by transfinite induction:
$V_0 = \emptyset$, $V_{\alpha+1} = \mathcal{P}(V_\alpha)$, for an ordinal $\alpha$ and $V_\alpha = \bigcup_{\beta < \alpha} V_\beta$, for a limit ordinal $\alpha$. Then $V=\bigcup_{\alpha \in Ord} V_\alpha$ is a class, called the universe of sets.
Given a set $M$, $\mathrm{def}(M)$ consists of all subsets $X\in \cP(M)$ such that there is a formula $\varphi(x,y_1,\ldots, y_n)$ and elements $a_1, \ldots, a_t\in M$ so that $X=\{ a\in M : M\models \varphi[a,a_1,\ldots, a_n]\}$ (see \cite[VI.3]{Eklof}).
The class of constructible sets is defined by transfinite induction: $L_0 = \emptyset$, $L_{\alpha+1} = \mathrm{def}(L_\alpha)$, for an ordinal $\alpha$ and $L_\alpha = \bigcup_{\beta < \alpha} L_\beta$, for a limit ordinal $\alpha$. Then $L=\bigcup_{\alpha \in Ord} L_\alpha$ is called the class of constructible sets. The Axiom of Constructibility says that every set is constructible, i.e. $V=L$. The result we are looking for is the following classical result found in Eklof-Mekler's book:

\begin{theorem}[{\cite[Theorem VII.1.4]{Eklof}}]\label{Ekolf}
    Assume $V=L$. If $R$ is not left perfect, then for every regular cardinal $\kappa$ which is not weakly compact there exists a $<\kappa$-generated $\kappa-free$, strongly $\kappa$-free module which is not projective.
\end{theorem}
The terminology needs some explanation. For an uncountable regular cardinal $\kappa$ a module $M$ is $<\kappa$-generated $\kappa$-free if there exists a family $\cF$ of  $<\kappa$-generated free submodules of $M$ such that every subset $X\subseteq M$ with $|X|<\kappa$ is contained in some $F\in\cF$ and  $\cF$ is closed under unions of increasing chains of length less than $\kappa$ (see \cite[Definition IV.1.1]{Eklof}). A $<\kappa$-generated $\kappa$-free module is $\kappa$-projective in our sense. The notion of strongly $\kappa$-free module is not of interest to us here (see \cite[Definition IV.1.8]{Eklof}). Combining Lemma \ref{necessary} and \cite[Corollary 2.2]{yilmaz} we conclude the following Corollary.
\begin{corollary}
Assume $V=L$. The stable category of $R$ has a weak generator  if and only if $R$ is left perfect.
\end{corollary}
\begin{proof}
In case $R$ is left perfect, the existence of a weak generator in the stable category has been shown  in \cite[Corollary 2.2]{yilmaz}. 
    Suppose $R$ is not left perfect,  and assume that  $G$ is a weak generator in the stable category of $R$. Set  $\kappa = \max(\aleph_1, |G|^+)$. Then $\kappa$, as a successor cardinal, is regular and not weakly compact. By Theorem \ref{Ekolf}, there exists a $\kappa$-projective module $M$ that is not projective, which contradicts Lemma \ref{necessary}.
\end{proof}

\section{Large-Cardinal Hypothesis}
In contrast to the Axiom of Constructibility, which provides the set-theoretic hypothesis for our nonexistence result, we now consider certain large-cardinal axioms that imply the existence of weak generators, for which we need some more terminology (see \cite[II.3]{Eklof}).
A filter on an infinite set $I$ is a family $D$ of subsets of $I$ that contains $I$ but not the empty set, is  upward closed under inclusion, and closed under finite intersections.
An ultrafilter on $I$ is a filter $D$ on $I$ such that for every $X\subseteq I$, either $X\in D$ or $I\setminus X\in D$. Equivalently, a filter $D$ is an ultrafilter on $I$ if and only if $D$ is a maximal filter.
 
For an infinite cardinal $\kappa$, a filter $D$ on $I$ is called $\kappa$-complete if
 it is closed under intersections of fewer than $\kappa$ of its members. 
 A cardinal $\kappa$ is called $L_{\lambda\omega}$-compact if for every set $I$, every $\kappa$-complete filter on $I$ is contained in a  $\lambda$-complete ultrafilter on $I$ (see \cite[Definition II.2.17]{Eklof}).\footnote{Recall that $\omega = \aleph_0$ denotes the cardinality of the natural numbers,while $\omega_1=\aleph_1 = \aleph_0^+$ is its successor.}

The existence of such $L_{\lambda\omega}$-compact cardinals $\kappa$ is a large-cardinal hypothesis beyond ZFC, with consequences for the existence of certain algebraic structures, which	 we wish to exploit. 

Recall that over a principal ideal domain (PID) any projective module is free. Thus any  $\kappa$-projective module is a $\kappa$-free module. We use the following theorem from Eklof-Mekler's book.  

\begin{theorem}[{\cite[Theorem II.3.10]{Eklof} }]\label{large}
Let $\lambda$ be an uncountable cardinal and $\kappa$ an $L_{\lambda\omega}$-compact cardinal. 
Let $R$ be a PID of cardinality $<\lambda$. Then any $\kappa$-free module is free.
\end{theorem}

Using Theorem \ref{large} and Lemma \ref{sufficient}, we can prove the existence of a weak generators in the stable category over a PID under the assumption of the existence of an $L_{\lambda\omega}$-compact cardinal.
\begin{corollary}
Let $R$ be a PID. If there exists a $L_{\lambda\omega}$-compact cardinal $\kappa$, for an uncountable cardinal $\lambda>|R|$, then the stable category of $R$ has a weak generator.
\end{corollary}

In particular, there exists a weak generator in the stable category of Abelian groups, if there exists a $L_{\omega_1\omega}$-compact cardinal $\kappa$, in which case one can take the direct sum of isomorphism classes of nonfree Abelian group of cardinality less than $\kappa$ as a weak generator.
Existence and properties of $L_{\omega_1\omega}$-compact cardinals $\kappa$ are discussed in \cite[Section 4.1]{BagariaMagidor}, where they are called $\omega_1$-strongly compact.

\medskip

As pointed out by Yilmaz Dur\u{g}un in a private conversation the problem of finding an $p$-indigent module is similar to looking for  Trlifaj's Whitehead test-modules. A module $G$ is called a Whitehead test-module for projetivity if whenever $\mathrm{Ext}^1_R(X,G)=0$, then $X$ is projective. In other words one is interested to show that $\Ker(\mathrm{Ext}^1(-,G))=\Proj(R)$. Our search for an $p$-indigent module $G$, or a weak generator in the stable category, is to show 
$\Ker(\sHom(G,-))=\Proj(R)$. It would be interesting to investigate further parallels.

\end{document}